\documentclass[times,numbers,12pt]{elsarticle}

\makeatletter
\def\ps@pprintTitle{%
 \let\@oddhead\@empty
 \let\@evenhead\@empty
 \def\@oddfoot{\centerline{\thepage}}%
 \let\@evenfoot\@oddfoot}
\makeatother

\usepackage{amsmath,amsfonts,amsthm,mathtools}
\usepackage[colorlinks, breaklinks=true]{hyperref}

\providecommand{\doi}[1]{\href{https://doi.org/#1}{doi:#1}}
\usepackage{xurl} 
\renewcommand{\doi}[1]{%
 \href{https://doi.org/#1}{\nolinkurl{doi:#1}}%
}

\usepackage{dsfont} 
\usepackage{xcolor} 
\usepackage{graphicx}
\usepackage{float}
\usepackage{geometry}
\theoremstyle{plain}
\newtheorem{theorem}{Theorem}

\newtheorem{remark}{Remark}

\newcommand{\N}{\mathbb{N}}
\newcommand{\R}{\mathbb{R}}

\newcommand{\PP}{\mathsf{P}} 
\newcommand{\bb}[1]{\boldsymbol{#1}}
\newcommand{\rd}{\hskip1pt\mathrm{d}\hskip0.75pt}
\newcommand{\ind}{\mathds{1}}
\newcommand{\OO}{\mathcal O}

\newcommand{\e}{\varepsilon}

\begin{document}

\begin{frontmatter}

\title{Normal integral representation for the joint survival function of the cumulative sums of the components of multinomial random vectors}

\author[a1]{Fr\'ed\'eric Ouimet}\ead{frederic.ouimet2@uqtr.ca}

\address[a1]{Universit\'e du Qu\'ebec \`a Trois-Rivi\`eres, Trois-Rivi\`eres, QC G8Z 4M3, Canada\vspace{-5mm}}

\begin{abstract}
This paper presents a multivariate normal integral representation for the joint survival function of the cumulative sums of the components of any multinomial random vector at interior lattice points. This result can be viewed as a multivariate analog of Equation~(7) in \citet{MR2154001}, whose proof starts from the beta integral representation of binomial survival probabilities and uses Laplace's method to improve Tusn\'ady's inequality. Our findings are based on a crucial relationship between the joint survival function of the cumulative sums of the components of any multinomial random vector and a Dirichlet probability over a corresponding cumulative-sum region. The main motivation is that such an explicit formula may eventually help streamline the conditional quantile-transformation arguments used in the multivariate KMT approximation of \citet{MR996984}, a connection left for future work. We provide numerical checks of the identity for $d = 2,3,4,5$.
\end{abstract}


\end{frontmatter}

\thispagestyle{empty}

\section{Introduction}\label{sec:intro}

One-dimensional quantile couplings between discrete distributions and Gaussian distributions are a central tool in probability and statistics. A classical example appears in the Koml\'os--Major--Tusn\'ady (KMT) approximation \citep[see, e.g.,][]{MR0375412,MR0402883,MR0666546,MR0972783,Major2000}, where the empirical distribution function is coupled with a Brownian bridge. As explained by \citet{MR2154001}, a key ingredient in that construction is the quantile coupling between a $\mathrm{Binomial}\hspace{0.2mm}(n,1/2)$ random variable and a normal random variable with the same mean and variance. More precisely, if $X\sim \mathrm{Binomial}\hspace{0.2mm}(n,1/2)$ and $Y\sim \mathcal{N}(n/2,n/4)$, then one defines cutpoints $-\infty = \beta_0 < \beta_1 < \cdots < \beta_n < \beta_{n+1} = \infty$ by
\[
\PP(X\geq k) = \PP(Y > \beta_k),\qquad k = 0,1,\ldots,n.
\]
Sharp control of the difference between the binomial quantiles and their Gaussian counterparts is the content of Tusn\'ady-type inequalities. The original inequality of \citet{Tusnady1977phd} and its later refinements \cite{MR1955348,MR2154001,MR4340237} show that the normal approximation to these quantiles is much more accurate than what follows from a crude central limit theorem.

The starting point of \citet{MR2154001} is not a local approximation of individual binomial probabilities, but rather an exact integral representation for the whole binomial tail:
\begin{equation}\label{eq:beta.integral.rep}
\PP(X\geq k) = \frac{n!}{(k-1)!(n-k)!}\int_0^{1/2} t^{k-1}(1-t)^{n-k}\rd t, \qquad 1 \leq k \leq n.
\end{equation}
This beta integral representation converts the discrete tail probability into a continuous integral. It is the one-dimensional prototype for the Dirichlet integral representation developed in the present paper.

Now, let us see how the Gaussian structure enters. For $n/2 < k < n$, let $K = k-1$, $N = n-1$, and $\e = (2K-N)/N$. Define
\[
2H(t) = (1+\e)\ln t+(1-\e)\ln(1-t),\qquad h(s) = H\left(\frac{1-s}{2}\right)-H(1/2),
\]
and let $\lambda_m$ denote the remainder term in Stirling's formula, as in \eqref{eq:log.Stirling} below. With
\[
\gamma(\e) = \frac{(1+\e)\ln(1+\e)+(1-\e)\ln(1-\e)-\e^2}{2\e^4}, \qquad \e \neq 0,
\]
where $\gamma(0) = 1/12$ by continuous extension, and
\[
a_N(\e) = \ln(1+N^{-1})+\lambda_N-\lambda_K-\lambda_{N-K}-\frac{1}{2}\ln(1-\e^2)-N\e^4\gamma(\e),
\]
Equation~(7) of \citet{MR2154001} can be written as
\[
\PP(X\geq k) = e^{a_N(\e)}\sqrt{\frac{N}{2\pi}}\int_0^1 \exp\left\{Nh(s)-\frac{N\e^2}{2}\right\}\rd s.
\]
Near its maximum at $s = 0$, the function $h$ starts as $h(s) = -\e s-s^2/2$ plus higher order terms. Thus the leading part of the integral is governed by the Gaussian expression
\[
\sqrt{\frac{N}{2\pi}}\int_0^\infty \exp\left\{-\frac{N(s+\e)^2}{2}\right\}\rd s = \overline{\Phi}(\e\sqrt{N}),
\]
where $\overline{\Phi}(x) = \PP(Z \geq x)$ for $Z\sim \mathcal{N}(0,1)$. \citet{MR2154001} then use this representation, together with careful Taylor bounds and inequalities for ratios of normal tail probabilities, to sharpen Tusn\'ady's inequality.

The goal of the present paper is to develop the corresponding integral representation in the multinomial setting, with the hope that such an explicit formula may eventually help simplify the difficult conditional quantile-transformation arguments used in the multivariate KMT approximation of \citet{MR996984}; see also \citet{MR1616527}. If $\bb{X} = (X_1,\ldots,X_d)\sim \mathrm{Multinomial}\hspace{0.2mm}(n,\bb{p})$, the natural multivariate analog of a binomial tail is the joint survival function of the cumulative sums of the components
\[
\PP(X_1+\cdots+X_i\geq k_1+\cdots+k_i, ~\forall i\in[d]),
\]
where $[d] \equiv \{1,\ldots,d\}$. These probabilities arise naturally when multinomial counts are constructed from a partition of the unit interval. In that construction, the event above can be expressed in terms of selected order statistics of independent uniform random variables; see \eqref{eq:order.statistics.link} below. The joint density of those order statistics yields a Dirichlet-type integral over a region $\mathcal{R}_d$ inside the simplex; see \eqref{eq:multinomial.dirichlet.relation} below. When $d = 1$ and $p_1 = 1/2$, this identity reduces exactly to the beta integral representation for binomial tails in \eqref{eq:beta.integral.rep}.

The main result of the paper (Theorem~\ref{thm:normal.integral.representation}) is an exact multivariate normal integral representation for this joint survival function at interior lattice points. This identity is not an asymptotic approximation. It is an exact rewriting of a Dirichlet integral, obtained by applying Stirling's formula and completing the square in the natural covariance structure of the multinomial distribution. In this sense, it is a multivariate analog of the normal integral representation used by \citet{MR2154001} in the binomial case to improve Tusn\'ady's inequality.

The rest of the paper is organized as follows. Section~\ref{sec:definitions.notation} introduces the necessary notation and derives the Dirichlet representation of the joint survival function of the cumulative sums of multinomial components. Section~\ref{sec:normal.integral.representation} states the multivariate normal integral representation. Section~\ref{sec:proofs} gives two proofs: the first expands the logarithm of the Dirichlet density directly, while the second uses a Laplace-type decomposition of the Dirichlet integral. Section~\ref{sec:numerical.verification} provides numerical checks of the claimed identity.

\section{Definitions and notation}\label{sec:definitions.notation}

For any integer $d\in \N$, the $d$-dimensional simplex and its interior are defined by
\[
\mathcal{S}_d = \big\{\bb{s}\in [0,1]^d: \|\bb{s}\|_1 \leq 1\big\}, \qquad \mathrm{Int}(\mathcal{S}_d) = \big\{\bb{s}\in (0,1)^d: \|\bb{s}\|_1 < 1\big\},
\]
where $\|\bb{s}\|_1 = \sum_{i=1}^d |s_i|$ denotes the $\ell_1$-norm in $\R^d$. Given a set of probability weights $\bb{p}\in \mathrm{Int}(\mathcal{S}_d)$, the $\mathrm{Multinomial}\hspace{0.2mm}(n,\bb{p})$ probability mass function is defined, for all $\bb{k}\in \N_0^d \cap n \mathcal{S}_d$, by
\[
p_n(\bb{k}) = \frac{n!}{(n - \|\bb{k}\|_1)! \prod_{i=1}^d k_i!} \, p_{d+1}^{n - \|\bb{k}\|_1} \prod_{i=1}^d p_i^{k_i},
\]
where $p_{d+1} = 1 - \|\bb{p}\|_1\in (0,1)$ and $n\in \N$. Throughout, when $\bb{s}\in \mathcal{S}_d$, set $s_{d+1} = 1 - \|\bb{s}\|_1$. The covariance matrix of the multinomial distribution is well-known to be $n \, \Sigma_{\bb{p}}$, where $\Sigma_{\bb{p}} = \text{diag}(\bb{p}) - \bb{p} \bb{p}^{\top}$; see, e.g., \cite[p.~377]{MR2168237}. From Theorem 1 of \citet{MR1157720}, it is also well known that $\det(\Sigma_{\bb{p}}) = p_1 \dots p_d \, p_{d+1}$. The centered multivariate normal density with covariance matrix $\Sigma_{\bb{p}}$ (the per-trial multinomial covariance matrix) is defined, for all $\bb{x}\in \R^d$, by
\[
\phi_{\Sigma_{\bb{p}}}(\bb{x}) = \frac{\exp\big(-\bb{x}^{\top} \Sigma_{\bb{p}}^{-1} \, \bb{x} / 2\big)}{\sqrt{(2\pi)^d \, \det(\Sigma_{\bb{p}})}}.
\]

Let $I_1 = (0,p_1]$ and $I_j = (p_1 + \dots + p_{j-1}, p_1 + \dots + p_j]$ for all $j\in [d]\backslash \{1\}$. If $U_1, \dots,U_n \smash{\stackrel{\mathrm{iid}}{\sim}} \mathrm{Uniform}(0,1)$, then
\[
\bb{X} = (X_1,\ldots,X_d) \equiv \Big(\sum_{i=1}^n \ind\{U_i\in I_1\},\ldots,\sum_{i=1}^n \ind\{U_i\in I_d\}\Big)\sim \mathrm{Multinomial}\hspace{0.2mm}(n,\bb{p}).
\]
For any $\bb{k}\in \N^d \cap n \mathcal{S}_d$, set $K_0 = 0$, $K_i = k_1 + \dots + k_i$ for all $i\in [d]$, and $K_{d+1} = n+1$. If $U_{(1)} \leq \dots \leq U_{(n)}$ denote the order statistics, then
\begin{equation}\label{eq:order.statistics.link}
\begin{aligned}
&\PP(X_1 + \dots + X_i \geq k_1 + \dots + k_i, ~\forall i\in [d]) \\[1mm]
&\qquad= \PP(U_{(K_i)} \leq p_1 + \dots + p_i, ~\forall i\in [d]) \\
&\qquad= \int_0^{p_1} \int_{u_1}^{p_1+p_2} \dots \int_{u_{d-1}}^{p_1+\dots+p_d} n! \prod_{i=1}^{d+1} \frac{(u_i - u_{i-1})^{K_i-K_{i-1}-1}}{(K_i - K_{i-1} - 1)!} \rd u_1 \rd u_2 \dots \rd u_d,
\end{aligned}
\end{equation}
where one defines $u_0 = 0$ and $u_{d+1} = 1$ in the multidimensional integral. After applying the change of variables $s_i = u_i - u_{i-1}$ for all $i\in [d]$, setting $s_{d+1} = 1 - \|\bb{s}\|_1$, and writing $k_{d+1} = K_{d+1} - K_d = (n + 1) - \|\bb{k}\|_1$, the above can be rewritten as
\[
\begin{aligned}
&\PP(X_1 + \dots + X_i \geq k_1 + \dots + k_i, ~\forall i\in [d]) \\
&\qquad= \int_0^{p_1} \int_0^{(p_1 - s_1) + p_2} \dots \int_0^{\sum_{k=1}^{d-1} (p_k - s_k) + p_d} n! \prod_{i=1}^{d+1} \frac{s_i^{k_i-1}}{(k_i - 1)!} \rd \bb{s}.
\end{aligned}
\]
This identity provides a direct relationship between the joint survival function of the cumulative sums of the components of any multinomial random vector and a Dirichlet probability over the corresponding cumulative-sum region.

In turn, letting $N = n - d$, $J_i = k_i - 1$ for all $i\in [d+1]$, so that $J_{d+1} = n - \|\bb{k}\|_1$, writing $\bb{J} = (J_1,\ldots,J_d)^{\top}$, and defining the regions
\begin{equation}\label{eq:def.R.d}
\mathcal{R}_d = \left\{\bb{s}\in \mathcal{S}_d : (s_1,s_2,\dots,s_i)\in \left(\sum_{k=1}^i p_k\right)\mathcal{S}_i ~~\forall i\in [d]\right\}, \qquad \mathcal{R}_d^{\circ} = \mathcal{R}_d \cap \mathrm{Int}(\mathcal{S}_d),
\end{equation}
one can write
\begin{equation}\label{eq:multinomial.dirichlet.relation}
\PP(X_1 + \dots + X_i \geq k_1 + \dots + k_i, ~\forall i\in [d])
= \int_{\mathcal{R}_d} \frac{(N + d)!}{N!} \times \frac{N!}{\prod_{i=1}^{d+1} J_i!} \prod_{i=1}^{d+1} s_i^{J_i} \rd \bb{s}.
\end{equation}
Here $\mathcal{R}_d^{\circ}$ denotes the part of $\mathcal{R}_d$ lying in $\mathrm{Int}(\mathcal{S}_d)$; it is not meant to denote the topological interior of $\mathcal{R}_d$, since some cumulative-sum constraints may be active. Since $\mathcal{R}_d\backslash \mathcal{R}_d^{\circ}$ has Lebesgue measure zero, either set may be used in \eqref{eq:multinomial.dirichlet.relation}; below we use $\mathcal{R}_d^{\circ}$ whenever logarithms of the $s_i$'s appear.

\section{Normal integral representation}\label{sec:normal.integral.representation}

For any integer $m\in \N$, let $\lambda_m$ denote the error term in Stirling's approximation for $\ln (m!)$:
\begin{equation}\label{eq:log.Stirling}
\ln(m!) = \frac{1}{2} \ln(2\pi m) + m \ln m - m + \lambda_m,
\end{equation}
where $(12 m + 1)^{-1} \leq \lambda_m \leq (12 m)^{-1}$; see, e.g., \citet{MR1502544}.

For the normal integral representation, assume from now on that
\[
\bb{k}\in \mathcal{K}_{n,d} = \left\{\bb{k}\in \N^d \cap n \mathcal{S}_d : k_i \geq 2 ~~\forall i\in [d],~\|\bb{k}\|_1 \leq n - 1\right\}.
\]
This condition is equivalent to $J_i\in \N$ for all $i\in [d+1]$, and it ensures that $N\in \N$ and that all logarithms below are finite. Together with the notation introduced in Section~\ref{sec:intro}, define, for all $i\in [d+1]$,
\[
\e_i = \frac{J_i / N - p_i}{p_i}, \qquad \widetilde{\e}_i = p_i \e_i = J_i/N - p_i,
\]
and write
\[
\bb{\e} = (\e_1,\ldots,\e_{d+1}), \qquad \widetilde{\bb{\e}} = (\widetilde{\e}_1,\ldots,\widetilde{\e}_d)^{\top}.
\]
Also, set
\[
\Lambda_N = \lambda_N - \sum_{i=1}^{d+1} \lambda_{J_i},
\]
and
\begin{align}
\Delta_N
&= \ln \left\{\frac{(N + d)!}{N! N^d}\right\} + \Lambda_N - \frac{1}{2} \sum_{i=1}^{d+1} \ln (1 + \e_i) - N \widetilde{\gamma}(\bb{\e}), \label{eq:Delta} \\
\widetilde{\gamma}(\bb{\e})
&= \sum_{i=1}^{d+1} p_i (1 + \e_i) \ln (1 + \e_i) - \frac{1}{2} \widetilde{\bb{\e}}^{\top} \Sigma_{\bb{p}}^{-1} \widetilde{\bb{\e}}, \label{eq:gamma.tilde} \\
\gamma^{\star}(\bb{s})
&= \sum_{i=1}^{d+1} p_i (1 + \e_i) \ln \left(\frac{s_i}{p_i}\right) - \left\{\widetilde{\bb{\e}}^{\top} \Sigma_{\bb{p}}^{-1} (\bb{s} - \bb{p}) - \frac{1}{2} (\bb{s} - \bb{p})^{\top} \Sigma_{\bb{p}}^{-1} (\bb{s} - \bb{p})\right\}. \label{eq:gamma.star}
\end{align}

Theorem~\ref{thm:normal.integral.representation} below expresses the joint survival function of the cumulative sums of the components of any multinomial random vector at interior lattice points $\bb{k}\in \mathcal{K}_{n,d}$ in terms of a multivariate normal integral. The result can be seen as a multivariate analog of Eq.~(7) of \citet{MR2154001}, who improved on Tusn\'ady's inequality. For an overview of the literature on Tusn\'ady's inequality and the most recent improvement in the bulk, see \citet{MR1955348} and \citet{MR4340237}, respectively.

\begin{theorem}[Normal integral representation]\label{thm:normal.integral.representation}
For all $\bb{k}\in \mathcal{K}_{n,d}$, one has
\[
\PP(X_1 + \dots + X_i \geq k_1 + \dots + k_i, ~\forall i\in [d]) = e^{\Delta_N} \int_{\mathcal{R}_d^{\circ}} \exp\left\{N \gamma^{\star}(\bb{s})\right\} N^{d/2} \phi_{\Sigma_{\bb{p}}}\{N^{1/2} (\bb{p} - \bb{s} + \widetilde{\bb{\e}})\} \rd \bb{s}.
\]
When $\max_{i\in [d+1]} |\e_i| \leq \eta < 1$, a useful expansion for $\widetilde{\gamma}(\bb{\e})$ can be found in \eqref{eq:gamma.epsilon.alternative}. An exact alternative expression for $\gamma^{\star}(\bb{s})$ can be found in \eqref{eq:gamma.star.alternative}.
\end{theorem}

\begin{remark}
The restriction $\bb{k}\in \mathcal{K}_{n,d}$ is used only for the normal representation above. The Dirichlet identity \eqref{eq:multinomial.dirichlet.relation} remains valid when some $J_i=0$, with the convention $0! = 1$, but those boundary cases are not covered by Theorem~\ref{thm:normal.integral.representation} because the logarithms in \eqref{eq:Delta}, \eqref{eq:gamma.tilde}, and \eqref{eq:gamma.star} would not all be finite.
\end{remark}

\begin{remark}
Theorem~\ref{thm:normal.integral.representation} also complements the local limit theorem for the multinomial distribution developed independently by \citet{MR750392} and \citet{MR4249129}.
\end{remark}

\section{Proofs}\label{sec:proofs}

\begin{proof}[\bf First proof of Theorem~\ref{thm:normal.integral.representation}]
Under the assumption $\bb{k}\in \mathcal{K}_{n,d}$, one has $J_i\in \N$ for all $i\in [d+1]$ and $\sum_{i=1}^{d+1} J_i = N$. For $\bb{s}\in \mathcal{R}_d^{\circ}$, by taking the logarithm of the integrand in \eqref{eq:multinomial.dirichlet.relation} and applying Stirling's formula \eqref{eq:log.Stirling}, one obtains
\begin{equation}\label{eq:log.p.before}
\begin{aligned}
&\ln \left\{\frac{(N + d)!}{N!} \times \frac{N!}{\prod_{i=1}^{d+1} J_i!} \prod_{i=1}^{d+1} s_i^{J_i}\right\} \\
&\quad= \ln \left\{\frac{(N + d)!}{N!}\right\} + \ln (N!) - \sum_{i=1}^{d+1} \ln (J_i!) + \sum_{i=1}^{d+1} J_i \ln s_i \\
&\quad= \ln \left\{\frac{(N + d)!}{N! N^d}\right\} - \frac{1}{2} \ln \left\{\left(\frac{2\pi}{N}\right)^d \, \prod_{i=1}^{d+1} (J_i/N)\right\} + \sum_{i=1}^{d+1} J_i \ln \left(\frac{s_i}{J_i / N}\right) + \lambda_N - \sum_{i=1}^{d+1} \lambda_{J_i} \\
&\quad= \ln \left\{\frac{(N + d)!}{N! N^d}\right\} - \frac{1}{2} \sum_{i=1}^{d+1} \ln (1 + \e_i) - \frac{1}{2} \ln \left\{\left(\frac{2\pi}{N}\right)^d \, \prod_{i=1}^{d+1} p_i\right\} \\
&\qquad+ N \sum_{i=1}^{d+1} p_i (1 + \e_i) \ln \left(\frac{s_i}{p_i}\right) + N \sum_{i=1}^{d+1} (J_i / N) \ln \left(\frac{p_i}{J_i / N}\right) + \Lambda_N.
\end{aligned}
\end{equation}
By the definition of $\widetilde{\gamma}(\bb{\e})$, one can write exactly
\begin{equation}\label{eq:entropy.1}
\sum_{i=1}^{d+1} (J_i / N) \ln \left(\frac{p_i}{J_i / N}\right) = - \frac{1}{2} \widetilde{\bb{\e}}^{\top} \Sigma_{\bb{p}}^{-1} \widetilde{\bb{\e}} - \widetilde{\gamma}(\bb{\e}).
\end{equation}
If, in addition, $\max_{i\in [d+1]} |\e_i| \leq \eta < 1$, then the Taylor expansion
\[
(1 + x) \ln(1 + x) = x + \frac{x^2}{2} - \frac{x^3}{6} + \frac{x^4}{12} + \OO_{\eta}(x^5), \quad |x| \leq \eta < 1,
\]
and the fact that $\widetilde{\e}_{d+1} = -\sum_{i=1}^d \widetilde{\e}_i$ show that $\widetilde{\gamma}(\bb{\e})$ can be expanded as follows:
\[
\begin{aligned}
\widetilde{\gamma}(\bb{\e})
&\equiv \sum_{i=1}^{d+1} p_i (1 + \e_i) \ln (1 + \e_i) - \frac{1}{2} \widetilde{\bb{\e}}^{\top} \Sigma_{\bb{p}}^{-1} \widetilde{\bb{\e}} \\
&= \frac{1}{2} \sum_{i,j=1}^d \widetilde{\e}_i \widetilde{\e}_j \left\{\frac{1}{p_i} \ind\{i = j\} + \frac{1}{p_{d+1}}\right\} - \frac{1}{2} \widetilde{\bb{\e}}^{\top} \Sigma_{\bb{p}}^{-1} \widetilde{\bb{\e}} - \frac{1}{6} \sum_{i,j,k=1}^d \widetilde{\e}_i \widetilde{\e}_j \widetilde{\e}_k \left\{\frac{1}{p_i^2} \ind\{i = j = k\} - \frac{1}{p_{d+1}^2}\right\} \\
&\quad+ \frac{1}{12} \sum_{i,j,k,\ell=1}^d \widetilde{\e}_i \widetilde{\e}_j \widetilde{\e}_k \widetilde{\e}_{\ell} \left\{\frac{1}{p_i^3} \ind\{i = j = k = \ell\} + \frac{1}{p_{d+1}^3}\right\} + \OO_{d,\bb{p},\eta}\big(\|\bb{\e}\|_1^5\big).
\end{aligned}
\]
Moreover, by \cite[Eq.~21]{MR1157720}, it is known that $(\Sigma_{\bb{p}}^{-1})_{ij} = p_i^{-1} \ind\{i = j\} + p_{d+1}^{-1}$ for all $i,j\in [d]$, so
\begin{equation}\label{eq:gamma.tilde.simplified}
\frac{1}{2} \sum_{i,j=1}^d \widetilde{\e}_i \widetilde{\e}_j \left\{\frac{1}{p_i} \ind\{i = j\} + \frac{1}{p_{d+1}}\right\} - \frac{1}{2} \widetilde{\bb{\e}}^{\top} \Sigma_{\bb{p}}^{-1} \widetilde{\bb{\e}} = 0.
\end{equation}
Under this additional condition, using the last two equations, one can rewrite $\widetilde{\gamma}(\bb{\e})$ as follows:
\begin{equation}\label{eq:gamma.epsilon.alternative}
\begin{aligned}
\widetilde{\gamma}(\bb{\e})
&=- \frac{1}{6} \sum_{i,j,k=1}^d \widetilde{\e}_i \widetilde{\e}_j \widetilde{\e}_k \left\{\frac{1}{p_i^2} \ind\{i = j = k\} - \frac{1}{p_{d+1}^2}\right\} \\
&\quad+ \frac{1}{12} \sum_{i,j,k,\ell=1}^d \widetilde{\e}_i \widetilde{\e}_j \widetilde{\e}_k \widetilde{\e}_{\ell} \left\{\frac{1}{p_i^3} \ind\{i = j = k = \ell\} + \frac{1}{p_{d+1}^3}\right\} + \OO_{d,\bb{p},\eta}\big(\|\bb{\e}\|_1^5\big).
\end{aligned}
\end{equation}
Now, consider $\bb{s} = (s_1,\ldots,s_d)\in \mathcal{R}_d^{\circ}$ as defined in \eqref{eq:def.R.d}. For all $i\in [d+1]$ and $t\in (0,1)$, let
\[
\delta_i(t) = \frac{t - p_i}{p_i}.
\]
In particular, note that this definition yields $\sum_{i=1}^{d+1} p_i \delta_i(s_i) = 0$. Upon applying Taylor's formula with Lagrange remainder to the function $x\mapsto \ln(1+x)$, there exist points $s_i^{\star}$, with $s_i^{\star}$ between $s_i$ and $p_i$, $i\in [d+1]$, such that
\[
\begin{aligned}
&\sum_{i=1}^{d+1} p_i (1 + \e_i) \ln \left(\frac{s_i}{p_i}\right)
= \sum_{i=1}^{d+1} p_i (1 + \e_i) \ln \{1 + \delta_i(s_i)\} \\
&\quad= \sum_{i=1}^{d+1} p_i (1 + \e_i) \delta_i(s_i) - \frac{1}{2} \sum_{i=1}^{d+1} p_i (1 + \e_i) \{\delta_i(s_i)\}^2 + \frac{1}{3} \sum_{i=1}^{d+1} p_i (1 + \e_i) \frac{\{\delta_i(s_i)\}^3}{\{1 + \delta_i(s_i^{\star})\}^3} \\
&\quad= \sum_{i=1}^{d+1} \widetilde{\e}_i \delta_i(s_i) - \frac{1}{2} \sum_{i=1}^{d+1} p_i \{\delta_i(s_i)\}^2 - \frac{1}{2} \sum_{i=1}^{d+1} \widetilde{\e}_i \{\delta_i(s_i)\}^2 + \frac{1}{3} \sum_{i=1}^{d+1} p_i (1 + \e_i) \frac{\{\delta_i(s_i)\}^3}{\{1 + \delta_i(s_i^{\star})\}^3} \\
&\quad= \widetilde{\bb{\e}}^{\top} \Sigma_{\bb{p}}^{-1} (\bb{s} - \bb{p}) - \frac{1}{2} (\bb{s} - \bb{p})^{\top} \Sigma_{\bb{p}}^{-1} (\bb{s} - \bb{p}) - \frac{1}{2} \sum_{i=1}^{d+1} \widetilde{\e}_i \{\delta_i(s_i)\}^2 + \frac{1}{3} \sum_{i=1}^{d+1} p_i (1 + \e_i) \frac{\{\delta_i(s_i)\}^3}{\{1 + \delta_i(s_i^{\star})\}^3},
\end{aligned}
\]
where $\widetilde{\bb{\e}} = (\widetilde{\e}_1,\ldots,\widetilde{\e}_d)$. Using the last equation, the quantity $\gamma^{\star}(\bb{s})$, as defined in \eqref{eq:gamma.star}, can be rewritten as
\begin{equation}\label{eq:gamma.star.alternative}
\begin{aligned}
\gamma^{\star}(\bb{s})
&\equiv \sum_{i=1}^{d+1} p_i (1 + \e_i) \ln \left(\frac{s_i}{p_i}\right) - \left\{\widetilde{\bb{\e}}^{\top} \Sigma_{\bb{p}}^{-1} (\bb{s} - \bb{p}) - \frac{1}{2} (\bb{s} - \bb{p})^{\top} \Sigma_{\bb{p}}^{-1} (\bb{s} - \bb{p})\right\} \\
&= - \frac{1}{2} \sum_{i=1}^{d+1} \widetilde{\e}_i \{\delta_i(s_i)\}^2 + \frac{1}{3} \sum_{i=1}^{d+1} p_i (1 + \e_i) \frac{\{\delta_i(s_i)\}^3}{\{1 + \delta_i(s_i^{\star})\}^3}.
\end{aligned}
\end{equation}
It readily follows that
\begin{equation}\label{eq:entropy.2}
\begin{aligned}
\sum_{i=1}^{d+1} p_i (1 + \e_i) \ln \left(\frac{s_i}{p_i}\right)
&= \widetilde{\bb{\e}}^{\top} \Sigma_{\bb{p}}^{-1} (\bb{s} - \bb{p}) - \frac{1}{2} (\bb{s} - \bb{p})^{\top} \Sigma_{\bb{p}}^{-1} (\bb{s} - \bb{p}) + \gamma^{\star}(\bb{s}) \\
&= \frac{1}{2} \widetilde{\bb{\e}}^{\top} \Sigma_{\bb{p}}^{-1} \widetilde{\bb{\e}} - \frac{1}{2} (\bb{s} - \bb{J}/N)^{\top} \Sigma_{\bb{p}}^{-1} (\bb{s} - \bb{J}/N) + \gamma^{\star}(\bb{s}).
\end{aligned}
\end{equation}
Therefore, by putting \eqref{eq:entropy.1} and \eqref{eq:entropy.2} back into \eqref{eq:log.p.before}, and exponentiating, one gets
\[
\begin{aligned}
\frac{(N + d)!}{N!} \times \frac{N!}{\prod_{i=1}^{d+1} J_i!} \prod_{i=1}^{d+1} s_i^{J_i}
&= \exp\left\{\Delta_N + N \gamma^{\star}(\bb{s})\right\} \frac{\exp\big\{-N (\bb{s} - \bb{J}/N)^{\top} \Sigma_{\bb{p}}^{-1} (\bb{s} - \bb{J}/N) / 2\big\}}{\sqrt{(2\pi / N)^d \det(\Sigma_{\bb{p}})}} \\
&= \exp\left\{\Delta_N + N \gamma^{\star}(\bb{s})\right\} N^{d/2} \phi_{\Sigma_{\bb{p}}}\{N^{1/2} (\bb{J}/N - \bb{s})\} \\[2mm]
&= \exp\left\{\Delta_N + N \gamma^{\star}(\bb{s})\right\} N^{d/2} \phi_{\Sigma_{\bb{p}}}\{N^{1/2} (\bb{p} - \bb{s} + \widetilde{\bb{\e}})\},
\end{aligned}
\]
where we recall that $\Delta_N = \ln \{(N + d)!/(N! N^d)\} + \Lambda_N - (1/2) \sum_{i=1}^{d+1} \ln (1 + \e_i) - N \widetilde{\gamma}(\bb{\e})$ from \eqref{eq:Delta}. Putting this last equation in \eqref{eq:multinomial.dirichlet.relation} yields
\[
\PP(X_1 + \dots + X_i \geq k_1 + \dots + k_i, ~\forall i\in [d])
= e^{\Delta_N} \int_{\mathcal{R}_d^{\circ}} \exp\left[N \gamma^{\star}(\bb{s})\right] N^{d/2} \phi_{\Sigma_{\bb{p}}}\{N^{1/2} (\bb{p} - \bb{s} + \widetilde{\bb{\e}})\} \rd \bb{s}.
\]
This concludes the proof.
\end{proof}

\begin{proof}[\bf Second proof of Theorem~\ref{thm:normal.integral.representation}]
For all $\bb{s}\in \mathrm{Int}(\mathcal{S}_d)$, with $s_{d+1} = 1 - \|\bb{s}\|_1$, define
\[
H(\bb{s}) = \sum_{i=1}^{d+1} p_i (1 + \e_i) \ln s_i = \sum_{i=1}^{d+1} \frac{J_i}{N} \ln s_i.
\]
For all $\bb{k}\in \mathcal{K}_{n,d}$, the integral representation of the joint survival function of the cumulative sums of the components of any multinomial random vector derived in \eqref{eq:multinomial.dirichlet.relation} can be rewritten as
\begin{equation}\label{eq:normal.integral.representation.eq.begin}
\PP(X_1 + \dots + X_i \geq k_1 + \dots + k_i, ~\forall i\in [d]) = \frac{(N + d)!}{N!} \times \frac{N!}{\prod_{i=1}^{d+1} J_i!} \int_{\mathcal{R}_d^{\circ}} \exp\left\{N H(\bb{s})\right\} \rd \bb{s}.
\end{equation}
By Stirling's formula \eqref{eq:log.Stirling}, note that
\[
m! = \sqrt{2\pi m} \, \exp(m \ln m - m + \lambda_m),
\]
where $(12 m + 1)^{-1} \leq \lambda_m \leq (12 m)^{-1}$. Since $\sum_{i=1}^{d+1} J_i = N$ and $\Lambda_N = \lambda_N - \sum_{i=1}^{d+1} \lambda_{J_i}$, one deduces that
\[
\frac{N!}{\prod_{i=1}^{d+1} J_i!} = \frac{\exp(\Lambda_N)}{\sqrt{(2\pi N)^d \prod_{i=1}^{d+1} \{p_i (1 + \e_i)\}}} \exp\left\{- N H(\bb{J}/N)\right\}.
\]
Therefore, one can rewrite \eqref{eq:normal.integral.representation.eq.begin} as
\begin{equation}\label{eq:normal.integral.representation.eq.next}
\begin{aligned}
&\PP(X_1 + \dots + X_i \geq k_1 + \dots + k_i, ~\forall i\in [d]) \\
&\qquad= \frac{(N + d)!}{N! N^d} \times \frac{N^{d/2} \exp(\Lambda_N)}{\sqrt{(2\pi)^d \prod_{i=1}^{d+1} \{p_i (1 + \e_i)\}}} \int_{\mathcal{R}_d^{\circ}} \exp\left[N \left\{H(\bb{s}) - H(\bb{J}/N)\right\}\right] \, \rd \bb{s}.
\end{aligned}
\end{equation}
Decompose the integral above as
\begin{equation}\label{eq:decomposition}
\int_{\mathcal{R}_d^{\circ}} \exp\left[N \left\{H(\bb{s}) - H(\bb{p})\right\} + N \left\{H(\bb{p}) - H(\bb{J}/N)\right\}\right] \, \rd \bb{s}.
\end{equation}
Using the quantity $\widetilde{\gamma}(\bb{\e})$ as defined in \eqref{eq:gamma.tilde}, one has
\begin{equation}\label{eq:max.diff.H}
H(\bb{p}) - H(\bb{J}/N)
= - \sum_{i=1}^{d+1} p_i (1 + \e_i) \ln(1 + \e_i)
= - \frac{1}{2} \widetilde{\bb{\e}}^{\top} \Sigma_{\bb{p}}^{-1} \widetilde{\bb{\e}} - \widetilde{\gamma}(\bb{\e}),
\end{equation}
where we recall that $\widetilde{\bb{\e}} = (p_1 \e_1,\ldots, p_d \e_d)^{\top}$, $\widetilde{\e}_{d+1} = p_{d+1} \e_{d+1}$, and $\Sigma_{\bb{p}}^{-1} = \{p_i^{-1} \ind\{i = j\} + p_{d+1}^{-1}\}_{1 \leq i,j \leq d}$.
Recall from \eqref{eq:Delta} that
\[
\Delta_N = \ln\left\{\frac{(N + d)!}{N! N^d}\right\} + \Lambda_N - \frac{1}{2} \sum_{i=1}^{d+1} \ln(1 + \e_i) - N \widetilde{\gamma}(\bb{\e}),
\]
so one can rewrite \eqref{eq:normal.integral.representation.eq.next} as
\[
\begin{aligned}
&\PP(X_1 + \dots + X_i \geq k_1 + \dots + k_i, ~\forall i\in [d]) \\[1mm]
&\qquad= \frac{N^{d/2} \exp(\Delta_N)}{\sqrt{(2\pi)^d \det(\Sigma_{\bb{p}})}} \int_{\mathcal{R}_d^{\circ}} \exp\left[N \left\{H(\bb{s}) - H(\bb{p})\right\} - \frac{N}{2} \widetilde{\bb{\e}}^{\top} \Sigma_{\bb{p}}^{-1} \widetilde{\bb{\e}}\right] \, \rd \bb{s}.
\end{aligned}
\]
Using \eqref{eq:gamma.star.alternative} and \eqref{eq:entropy.2}, one has
\[
\begin{aligned}
H(\bb{s}) - H(\bb{p})
&= \sum_{i=1}^{d+1} \frac{J_i}{N} \ln (s_i) - \sum_{i=1}^{d+1} \frac{J_i}{N} \ln (p_i)
= \sum_{i=1}^{d+1} p_i (1 + \e_i) \ln \left(\frac{s_i}{p_i}\right) \\
&= \widetilde{\bb{\e}}^{\top} \Sigma_{\bb{p}}^{-1} (\bb{s} - \bb{p}) - \frac{1}{2} (\bb{s} - \bb{p})^{\top} \Sigma_{\bb{p}}^{-1} (\bb{s} - \bb{p}) - \frac{1}{2} \sum_{i=1}^{d+1} \widetilde{\e}_i \{\delta_i(s_i)\}^2 + \frac{1}{3} \sum_{i=1}^{d+1} p_i (1 + \e_i) \frac{\{\delta_i(s_i)\}^3}{\{1 + \delta_i(s_i^{\star})\}^3} \\
&= \widetilde{\bb{\e}}^{\top} \Sigma_{\bb{p}}^{-1} (\bb{s} - \bb{p}) - \frac{1}{2} (\bb{s} - \bb{p})^{\top} \Sigma_{\bb{p}}^{-1} (\bb{s} - \bb{p}) + \gamma^{\star}(\bb{s}) \\
&= \frac{1}{2} \widetilde{\bb{\e}}^{\top} \Sigma_{\bb{p}}^{-1} \widetilde{\bb{\e}} - \frac{1}{2} (\bb{s} - \bb{J}/N)^{\top} \Sigma_{\bb{p}}^{-1} (\bb{s} - \bb{J}/N) + \gamma^{\star}(\bb{s}).
\end{aligned}
\]
Putting this in the previous equation yields
\[
\begin{aligned}
\PP(X_1 + \dots + X_i \geq k_1 + \dots + k_i, ~\forall i\in [d])
&= e^{\Delta_N} \int_{\mathcal{R}_d^{\circ}} \exp\left\{N \gamma^{\star}(\bb{s})\right\} N^{d/2} \phi_{\Sigma_{\bb{p}}}\{N^{1/2} (\bb{J}/N - \bb{s})\} \rd \bb{s} \\
&= e^{\Delta_N} \int_{\mathcal{R}_d^{\circ}} \exp\left\{N \gamma^{\star}(\bb{s})\right\} N^{d/2} \phi_{\Sigma_{\bb{p}}}\{N^{1/2} (\bb{p} - \bb{s} + \widetilde{\bb{\e}})\} \rd \bb{s}.
\end{aligned}
\]
This concludes the proof.
\end{proof}

\begin{remark}
The choice of decomposition in \eqref{eq:decomposition} can be explained as follows. The Hessian matrix of $H$ at $\bb{s}\in \mathrm{Int}(\mathcal{S}_d)$ is equal to
\[
\left(-\frac{J_i}{N s_i^2} \ind\{i = j\} - \frac{J_{d+1}}{N s_{d+1}^2}\right)_{1 \leq i,j \leq d},
\]
which is (symmetric) negative definite, so that $H$ is strictly concave on $\mathrm{Int}(\mathcal{S}_d)$. Given the strict concavity and the fact that $\overline{\bb{s}} = \bb{J}/N\in \mathrm{Int}(\mathcal{S}_d)$ is the unique point satisfying
\[
\nabla H(\overline{\bb{s}}) = \left(\frac{J_i}{N \, \overline{s}_i} - \frac{J_{d+1}}{N \, \overline{s}_{d+1}}\right)_{1 \leq i \leq d} = \bb{0}
\]
in the interior of the simplex, the global maximum of $H$ on $\mathrm{Int}(\mathcal{S}_d)$ is achieved at $\bb{J}/N$. If $\bb{J}/N\in \mathcal{R}_d^{\circ}$, then \eqref{eq:decomposition} is the usual centering at the unconstrained maximum in Laplace's method. If $\bb{J}/N\notin \mathcal{R}_d^{\circ}$, then the constrained maximum over $\mathcal{R}_d$ is attained on the boundary of $\mathcal{R}_d$ and need not be equal to $\bb{p}$; in that case, \eqref{eq:decomposition} should be viewed simply as an exact algebraic decomposition.
\end{remark}

\section{Numerical checks}\label{sec:numerical.verification}

We numerically check the exact identity in Theorem~\ref{thm:normal.integral.representation} by comparing, for randomly generated admissible triples $(n,\bb{p},\bb{k})$, the two quantities
\[
\begin{aligned}
L(n,\bb{p},\bb{k}) &= \PP(X_1 + \dots + X_i \geq k_1 + \dots + k_i, ~\forall i\in [d]), \\
R(n,\bb{p},\bb{k}) &= e^{\Delta_N} \int_{\mathcal{R}_d^{\circ}} \exp\left\{N\gamma^{\star}(\bb{s})\right\} N^{d/2} \phi_{\Sigma_{\bb{p}}} \{N^{1/2}(\bb{p}-\bb{s}+\widetilde{\bb{\e}})\} \rd\bb{s}.
\end{aligned}
\]
The first quantity, $L(n,\bb{p},\bb{k})$, is evaluated by direct summation of the multinomial probabilities,
\[
L(n,\bb{p},\bb{k}) = \sum_{\substack{x_1,\ldots,x_{d+1}\in\N_0 \\ x_1+\cdots+x_{d+1} = n}} \frac{n!}{x_1!\cdots x_{d+1}!} \prod_{j=1}^{d+1}p_j^{x_j} \ind\{x_1+\cdots+x_i \geq k_1+\cdots+k_i,~\forall i\in[d]\},
\]
where $p_{d+1} = 1-\|\bb{p}\|_1$. This is an exact finite sum, up to floating-point roundoff; no Monte Carlo approximation is used.

The second quantity is evaluated by deterministic adaptive numerical integration. To integrate over $\mathcal{R}_d$, the code uses cumulative coordinates. Let $P_i = p_1+\cdots+p_i$ for $i\in[d]$, and map $\bb{u}\in[0,1]^d$ to $\bb{s}\in\mathcal{R}_d$ recursively by
\vspace{-2mm}
\[
q_0 = 0,\qquad
q_i = q_{i-1}+u_i(P_i-q_{i-1}),\qquad
s_i = q_i-q_{i-1},\qquad i\in[d].
\]
The Jacobian of this transformation is $\prod_{i=1}^d (P_i-q_{i-1})$, and the resulting integral over $[0,1]^d$ is computed using adaptive cubature. The comparison therefore checks the equality in Theorem~\ref{thm:normal.integral.representation} by comparing an exact multinomial summation against a deterministic numerical integral.

For each tested case, we report the absolute and symmetric relative errors:
\[
\mathrm{AbsErr} = |R(n,\bb{p},\bb{k})-L(n,\bb{p},\bb{k})|, \qquad
\mathrm{RelErr} = \frac{2|R(n,\bb{p},\bb{k})-L(n,\bb{p},\bb{k})|}{|R(n,\bb{p},\bb{k})|+|L(n,\bb{p},\bb{k})|}.
\]
The symmetric relative error is used to avoid favoring either side of the identity.

The numerical experiment was run for $d = 2,3,4,5$; see Tables~\ref{tab:theorem-checks-d2}--\ref{tab:theorem-checks-d5}. For each dimension, $40$ randomized cases were generated. In each case, $n$ was sampled uniformly from a finite integer interval, and the full probability vector
\[
\bb{p}^{+} = (p_1,\ldots,p_d,p_{d+1})
\]
was sampled from a symmetric Dirichlet distribution with parameter $2$, conditional on all components being at least $0.03$. The vector $\bb{k}$ was sampled subject to the theorem's restrictions $k_i\geq 2$ for all $i\in[d]$ and $\|\bb{k}\|_1\leq n-1$. More precisely, after assigning the baseline value $2$ to each component of $\bb{k}$, the remaining slack $n-1-2d$ was randomly allocated among $d+1$ cells using probabilities sampled from a symmetric Dirichlet distribution with parameter $1$, and only the first $d$ allocations were added to $\bb{k}$. Thus every generated case satisfies $\bb{k}\in\mathcal{K}_{n,d}$. The ranges for $n$ and the cubature tolerances were set to $12\leq n\leq 35$ and $10^{-7}$, respectively, for all $d = 2,3,4,5$.

The \textsf{R} code that generated Tables~\ref{tab:theorem-checks-d2}--\ref{tab:theorem-checks-d5} is available online in the GitHub repository \cite{Ouimet2026github}.


\begin{table}[H]
\centering
\caption{Results from 40 randomized checks of the main theorem for $d=2$. Here $L$ denotes the left side and $R$ denotes the right side of the equation in Theorem~\ref{thm:normal.integral.representation}, and the relative error is $2|R-L|/(|R|+|L|)$.}
\label{tab:theorem-checks-d2}
\vspace{2em}
\begin{tabular}{r r l l r r}
\hline
Case & $n$ & $\bb{p}^+$ & $\bb{k}$ & Abs. error & Rel. error \\
\hline
1 & 33 & $(0.277, 0.209, 0.513)$ & $(7, 13)$ & 8.381e-12 & 7.330e-11 \\
2 & 14 & $(0.149, 0.575, 0.276)$ & $(3, 6)$ & 6.537e-12 & 2.063e-11 \\
3 & 16 & $(0.249, 0.221, 0.530)$ & $(3, 8)$ & 7.723e-12 & 1.173e-10 \\
4 & 28 & $(0.447, 0.240, 0.313)$ & $(13, 4)$ & 2.864e-11 & 5.854e-11 \\
5 & 25 & $(0.184, 0.664, 0.153)$ & $(9, 10)$ & 6.767e-12 & 2.426e-10 \\
6 & 33 & $(0.071, 0.478, 0.451)$ & $(22, 7)$ & 2.177e-19 & \textbf{1.961e-03} \\
7 & 33 & $(0.156, 0.614, 0.230)$ & $(4, 21)$ & 5.533e-10 & 1.023e-09 \\
8 & 17 & $(0.200, 0.414, 0.386)$ & $(2, 5)$ & 2.721e-11 & 3.140e-11 \\
9 & 14 & $(0.301, 0.293, 0.406)$ & $(2, 6)$ & 5.970e-10 & 8.992e-10 \\
10 & 35 & $(0.310, 0.266, 0.424)$ & $(15, 12)$ & 1.840e-12 & 2.551e-10 \\
11 & 13 & $(0.624, 0.159, 0.218)$ & $(4, 2)$ & 2.399e-11 & 2.414e-11 \\
12 & 31 & $(0.524, 0.083, 0.393)$ & $(18, 10)$ & 3.518e-13 & 1.287e-09 \\
13 & 25 & $(0.120, 0.451, 0.428)$ & $(2, 3)$ & 1.060e-10 & 1.294e-10 \\
14 & 32 & $(0.061, 0.741, 0.198)$ & $(13, 11)$ & 8.192e-18 & 4.330e-10 \\
15 & 31 & $(0.108, 0.343, 0.549)$ & $(4, 9)$ & 2.185e-10 & 6.059e-10 \\
16 & 25 & $(0.126, 0.340, 0.535)$ & $(4, 14)$ & 3.205e-11 & 4.621e-09 \\
17 & 31 & $(0.479, 0.127, 0.394)$ & $(12, 11)$ & 4.244e-10 & 5.094e-09 \\
18 & 12 & $(0.413, 0.249, 0.338)$ & $(4, 5)$ & 1.835e-12 & 5.120e-12 \\
19 & 30 & $(0.082, 0.467, 0.451)$ & $(16, 9)$ & 5.537e-16 & \textbf{1.097e-05} \\
20 & 19 & $(0.119, 0.753, 0.128)$ & $(2, 15)$ & 1.032e-10 & 2.596e-10 \\
21 & 23 & $(0.370, 0.303, 0.326)$ & $(4, 10)$ & 5.045e-10 & 6.227e-10 \\
22 & 21 & $(0.208, 0.490, 0.301)$ & $(14, 6)$ & 8.603e-17 & 5.945e-11 \\
23 & 28 & $(0.622, 0.095, 0.283)$ & $(11, 3)$ & 1.803e-12 & 1.815e-12 \\
24 & 19 & $(0.162, 0.596, 0.242)$ & $(7, 5)$ & 1.020e-11 & 4.208e-10 \\
25 & 32 & $(0.185, 0.386, 0.430)$ & $(11, 8)$ & 6.337e-11 & 2.984e-09 \\
26 & 14 & $(0.399, 0.172, 0.428)$ & $(2, 6)$ & 7.570e-12 & 1.238e-11 \\
27 & 29 & $(0.348, 0.340, 0.312)$ & $(10, 18)$ & 1.324e-13 & 5.022e-10 \\
28 & 12 & $(0.435, 0.427, 0.138)$ & $(4, 2)$ & 6.192e-11 & 7.354e-11 \\
29 & 33 & $(0.274, 0.558, 0.168)$ & $(3, 22)$ & 3.320e-10 & 3.658e-10 \\
30 & 18 & $(0.638, 0.153, 0.209)$ & $(3, 13)$ & 1.484e-10 & 6.111e-10 \\
31 & 15 & $(0.230, 0.357, 0.413)$ & $(3, 6)$ & 2.005e-11 & 4.305e-11 \\
32 & 18 & $(0.143, 0.485, 0.372)$ & $(8, 8)$ & 1.969e-14 & 7.514e-11 \\
33 & 28 & $(0.312, 0.524, 0.165)$ & $(22, 4)$ & 3.349e-17 & 1.200e-10 \\
34 & 35 & $(0.244, 0.334, 0.423)$ & $(12, 15)$ & 1.212e-11 & 1.734e-09 \\
35 & 19 & $(0.170, 0.455, 0.375)$ & $(11, 2)$ & 5.249e-15 & 7.843e-11 \\
36 & 29 & $(0.158, 0.605, 0.237)$ & $(12, 6)$ & 1.258e-12 & 1.412e-09 \\
37 & 14 & $(0.374, 0.268, 0.358)$ & $(3, 8)$ & 5.201e-11 & 2.590e-10 \\
38 & 23 & $(0.284, 0.244, 0.472)$ & $(3, 3)$ & 1.473e-10 & 1.510e-10 \\
39 & 12 & $(0.123, 0.500, 0.377)$ & $(3, 5)$ & 1.303e-11 & 1.042e-10 \\
40 & 25 & $(0.386, 0.423, 0.191)$ & $(2, 13)$ & 1.956e-09 & 1.964e-09 \\
\hline
\end{tabular}
\par\vspace{0.5em}
Errors larger than $10^{-6}$ are highlighted in bold.
\end{table}


\begin{table}[H]
\centering
\caption{Results from 40 randomized checks of the main theorem for $d=3$. Here $L$ denotes the left side and $R$ denotes the right side of the equation in Theorem~\ref{thm:normal.integral.representation}, and the relative error is $2|R-L|/(|R|+|L|)$.}
\label{tab:theorem-checks-d3}
\vspace{2em}
\begin{tabular}{r r l l r r}
\hline
Case & $n$ & $\bb{p}^+$ & $\bb{k}$ & Abs. error & Rel. error \\
\hline
1 & 35 & $(0.145, 0.071, 0.242, 0.542)$ & $(2, 15, 12)$ & 6.359e-16 & 5.696e-10 \\
2 & 16 & $(0.334, 0.173, 0.383, 0.110)$ & $(2, 3, 2)$ & 7.017e-11 & 7.320e-11 \\
3 & 35 & $(0.291, 0.185, 0.480, 0.044)$ & $(15, 2, 5)$ & 1.116e-11 & 1.957e-10 \\
4 & 21 & $(0.112, 0.555, 0.267, 0.065)$ & $(5, 5, 7)$ & 4.691e-12 & 6.016e-11 \\
5 & 32 & $(0.341, 0.151, 0.171, 0.337)$ & $(2, 12, 13)$ & 1.802e-11 & 9.159e-10 \\
6 & 30 & $(0.309, 0.266, 0.067, 0.357)$ & $(19, 5, 5)$ & 3.024e-15 & 2.121e-09 \\
7 & 25 & $(0.041, 0.470, 0.144, 0.345)$ & $(7, 7, 4)$ & 1.099e-14 & 3.342e-10 \\
8 & 13 & $(0.228, 0.255, 0.277, 0.240)$ & $(5, 4, 2)$ & 2.516e-12 & 5.962e-11 \\
9 & 28 & $(0.080, 0.278, 0.428, 0.214)$ & $(3, 11, 4)$ & 1.478e-11 & 2.506e-10 \\
10 & 21 & $(0.082, 0.123, 0.693, 0.103)$ & $(3, 6, 2)$ & 3.335e-12 & 2.571e-10 \\
11 & 22 & $(0.179, 0.073, 0.195, 0.553)$ & $(2, 9, 5)$ & 1.474e-13 & 7.438e-11 \\
12 & 23 & $(0.432, 0.267, 0.070, 0.230)$ & $(6, 6, 8)$ & 1.693e-10 & 8.968e-10 \\
13 & 34 & $(0.351, 0.366, 0.156, 0.127)$ & $(9, 4, 9)$ & 1.335e-10 & 1.494e-10 \\
14 & 18 & $(0.077, 0.437, 0.324, 0.163)$ & $(4, 5, 2)$ & 6.677e-13 & 1.671e-11 \\
15 & 33 & $(0.214, 0.503, 0.107, 0.176)$ & $(8, 2, 8)$ & 5.313e-12 & 1.291e-11 \\
16 & 18 & $(0.293, 0.302, 0.184, 0.221)$ & $(4, 3, 5)$ & 2.904e-11 & 3.783e-11 \\
17 & 17 & $(0.250, 0.208, 0.191, 0.351)$ & $(4, 5, 3)$ & 4.307e-11 & 1.822e-10 \\
18 & 18 & $(0.305, 0.064, 0.107, 0.525)$ & $(5, 7, 3)$ & 1.099e-12 & 9.402e-10 \\
19 & 14 & $(0.198, 0.189, 0.274, 0.339)$ & $(3, 6, 3)$ & 1.789e-11 & 8.474e-10 \\
20 & 29 & $(0.194, 0.364, 0.309, 0.133)$ & $(15, 4, 7)$ & 6.473e-15 & 8.165e-11 \\
21 & 31 & $(0.115, 0.146, 0.446, 0.294)$ & $(7, 6, 7)$ & 3.968e-13 & 2.582e-11 \\
22 & 20 & $(0.579, 0.200, 0.053, 0.169)$ & $(7, 4, 4)$ & 2.662e-10 & 2.996e-10 \\
23 & 23 & $(0.545, 0.204, 0.144, 0.107)$ & $(7, 5, 2)$ & 3.338e-10 & 3.369e-10 \\
24 & 31 & $(0.258, 0.461, 0.151, 0.130)$ & $(5, 5, 13)$ & 6.744e-11 & 7.339e-11 \\
25 & 17 & $(0.282, 0.167, 0.181, 0.370)$ & $(10, 2, 3)$ & 2.302e-13 & 1.425e-10 \\
26 & 34 & $(0.128, 0.221, 0.188, 0.463)$ & $(10, 3, 6)$ & 1.840e-12 & 2.468e-10 \\
27 & 25 & $(0.238, 0.375, 0.297, 0.089)$ & $(8, 6, 7)$ & 2.777e-11 & 1.320e-10 \\
28 & 15 & $(0.122, 0.443, 0.211, 0.224)$ & $(3, 8, 3)$ & 5.945e-12 & 2.187e-10 \\
29 & 12 & $(0.101, 0.237, 0.433, 0.228)$ & $(4, 3, 2)$ & 1.407e-12 & 1.369e-10 \\
30 & 19 & $(0.388, 0.173, 0.193, 0.246)$ & $(5, 2, 7)$ & 2.510e-10 & 3.858e-10 \\
31 & 22 & $(0.089, 0.113, 0.360, 0.438)$ & $(5, 9, 6)$ & 4.910e-17 & 5.126e-11 \\
32 & 29 & $(0.115, 0.365, 0.371, 0.149)$ & $(8, 9, 4)$ & 2.164e-12 & 2.750e-10 \\
33 & 14 & $(0.664, 0.102, 0.084, 0.150)$ & $(3, 6, 2)$ & 6.641e-11 & 8.021e-11 \\
34 & 35 & $(0.351, 0.111, 0.314, 0.223)$ & $(4, 10, 3)$ & 2.316e-11 & 2.829e-11 \\
35 & 22 & $(0.330, 0.215, 0.300, 0.156)$ & $(2, 4, 12)$ & 1.628e-12 & 2.177e-12 \\
36 & 20 & $(0.320, 0.252, 0.333, 0.095)$ & $(11, 3, 4)$ & 9.725e-13 & 4.943e-11 \\
37 & 16 & $(0.257, 0.244, 0.214, 0.285)$ & $(5, 6, 3)$ & 1.994e-11 & 4.932e-10 \\
38 & 35 & $(0.294, 0.344, 0.323, 0.039)$ & $(5, 11, 11)$ & 2.180e-11 & 2.223e-11 \\
39 & 25 & $(0.480, 0.269, 0.057, 0.195)$ & $(9, 8, 6)$ & 8.902e-11 & 8.243e-10 \\
40 & 13 & $(0.420, 0.194, 0.320, 0.066)$ & $(3, 4, 4)$ & 7.452e-11 & 9.735e-11 \\
\hline
\end{tabular}
\par\vspace{0.5em}
Errors larger than $10^{-6}$ are highlighted in bold.
\end{table}


\begin{table}[H]
\centering
\caption{Results from 40 randomized checks of the main theorem for $d=4$. Here $L$ denotes the left side and $R$ denotes the right side of the equation in Theorem~\ref{thm:normal.integral.representation}, and the relative error is $2|R-L|/(|R|+|L|)$.}
\label{tab:theorem-checks-d4}
\vspace{2em}
\begin{tabular}{r r l l r r}
\hline
Case & $n$ & $\bb{p}^+$ & $\bb{k}$ & Abs. error & Rel. error \\
\hline
1 & 25 & $(0.084, 0.184, 0.259, 0.406, 0.067)$ & $(4, 6, 3, 10)$ & 1.027e-12 & 2.539e-11 \\
2 & 33 & $(0.079, 0.323, 0.277, 0.186, 0.135)$ & $(2, 15, 13, 2)$ & 5.760e-14 & 6.661e-11 \\
3 & 17 & $(0.260, 0.191, 0.046, 0.292, 0.211)$ & $(2, 5, 6, 3)$ & 8.684e-13 & 7.474e-11 \\
4 & 18 & $(0.388, 0.248, 0.227, 0.054, 0.083)$ & $(5, 3, 2, 3)$ & 8.157e-10 & 9.312e-10 \\
5 & 23 & $(0.179, 0.233, 0.159, 0.297, 0.132)$ & $(2, 2, 12, 5)$ & 4.141e-12 & 3.765e-11 \\
6 & 25 & $(0.267, 0.074, 0.112, 0.280, 0.267)$ & $(3, 3, 2, 14)$ & 6.534e-11 & 9.657e-10 \\
7 & 12 & $(0.135, 0.127, 0.417, 0.057, 0.263)$ & $(2, 2, 2, 4)$ & 1.283e-11 & 8.436e-11 \\
8 & 14 & $(0.209, 0.137, 0.087, 0.476, 0.091)$ & $(3, 2, 6, 2)$ & 4.103e-13 & 5.670e-11 \\
9 & 16 & $(0.213, 0.367, 0.140, 0.211, 0.069)$ & $(2, 3, 2, 2)$ & 1.413e-10 & 1.608e-10 \\
10 & 29 & $(0.080, 0.260, 0.407, 0.223, 0.031)$ & $(3, 9, 2, 6)$ & 1.385e-11 & 8.685e-11 \\
11 & 34 & $(0.498, 0.196, 0.176, 0.079, 0.051)$ & $(6, 5, 11, 3)$ & 6.169e-09 & 6.170e-09 \\
12 & 27 & $(0.286, 0.235, 0.164, 0.198, 0.117)$ & $(3, 7, 13, 3)$ & 1.911e-12 & 7.615e-11 \\
13 & 27 & $(0.153, 0.186, 0.142, 0.139, 0.380)$ & $(9, 7, 8, 2)$ & 2.024e-16 & 2.726e-10 \\
14 & 18 & $(0.237, 0.193, 0.155, 0.209, 0.206)$ & $(2, 2, 4, 2)$ & 2.090e-10 & 2.363e-10 \\
15 & 18 & $(0.390, 0.188, 0.052, 0.207, 0.162)$ & $(4, 3, 4, 6)$ & 9.655e-12 & 5.789e-11 \\
16 & 31 & $(0.277, 0.045, 0.287, 0.316, 0.076)$ & $(7, 3, 7, 10)$ & 1.613e-10 & 3.290e-10 \\
17 & 17 & $(0.081, 0.079, 0.437, 0.128, 0.274)$ & $(6, 2, 5, 3)$ & 5.173e-15 & 6.530e-11 \\
18 & 26 & $(0.234, 0.285, 0.042, 0.279, 0.159)$ & $(3, 4, 8, 2)$ & 2.816e-10 & 5.472e-10 \\
19 & 20 & $(0.260, 0.291, 0.151, 0.107, 0.191)$ & $(2, 2, 5, 6)$ & 3.691e-10 & 4.500e-10 \\
20 & 18 & $(0.213, 0.041, 0.201, 0.227, 0.318)$ & $(3, 2, 6, 3)$ & 3.540e-12 & 4.186e-11 \\
21 & 22 & $(0.237, 0.339, 0.110, 0.246, 0.069)$ & $(4, 8, 3, 3)$ & 7.543e-12 & 1.545e-11 \\
22 & 16 & $(0.139, 0.227, 0.127, 0.187, 0.319)$ & $(7, 3, 2, 3)$ & 1.866e-14 & 7.125e-11 \\
23 & 28 & $(0.047, 0.188, 0.374, 0.321, 0.071)$ & $(3, 2, 10, 11)$ & 6.638e-12 & 6.884e-11 \\
24 & 16 & $(0.199, 0.236, 0.260, 0.192, 0.113)$ & $(2, 4, 2, 6)$ & 1.233e-11 & 2.237e-11 \\
25 & 21 & $(0.242, 0.110, 0.210, 0.385, 0.053)$ & $(5, 3, 2, 9)$ & 4.149e-11 & 1.092e-10 \\
26 & 19 & $(0.052, 0.110, 0.330, 0.238, 0.269)$ & $(3, 9, 4, 2)$ & 1.244e-16 & 1.745e-10 \\
27 & 18 & $(0.052, 0.232, 0.238, 0.339, 0.140)$ & $(3, 6, 5, 2)$ & 7.462e-14 & 3.526e-11 \\
28 & 30 & $(0.142, 0.149, 0.311, 0.144, 0.254)$ & $(7, 3, 2, 17)$ & 6.329e-14 & 1.229e-10 \\
29 & 23 & $(0.089, 0.276, 0.173, 0.290, 0.173)$ & $(2, 6, 4, 6)$ & 6.762e-12 & 1.935e-11 \\
30 & 25 & $(0.216, 0.128, 0.146, 0.166, 0.344)$ & $(3, 4, 8, 2)$ & 5.779e-11 & 3.415e-10 \\
31 & 32 & $(0.387, 0.031, 0.239, 0.064, 0.279)$ & $(2, 5, 11, 13)$ & 3.406e-13 & 8.799e-10 \\
32 & 31 & $(0.292, 0.158, 0.254, 0.256, 0.040)$ & $(3, 4, 8, 9)$ & 3.664e-09 & 3.690e-09 \\
33 & 27 & $(0.105, 0.146, 0.401, 0.158, 0.189)$ & $(3, 9, 3, 4)$ & 5.993e-12 & 2.849e-10 \\
34 & 29 & $(0.241, 0.237, 0.279, 0.164, 0.079)$ & $(9, 7, 9, 2)$ & 6.904e-12 & 1.598e-10 \\
35 & 14 & $(0.161, 0.213, 0.388, 0.159, 0.079)$ & $(3, 2, 2, 4)$ & 2.260e-12 & 6.552e-12 \\
36 & 34 & $(0.264, 0.100, 0.468, 0.075, 0.094)$ & $(2, 4, 4, 14)$ & 1.962e-09 & 1.972e-09 \\
37 & 13 & $(0.290, 0.232, 0.047, 0.267, 0.164)$ & $(2, 3, 3, 4)$ & 6.687e-12 & 2.799e-11 \\
38 & 30 & $(0.180, 0.468, 0.173, 0.137, 0.041)$ & $(4, 7, 10, 2)$ & 5.290e-10 & 6.651e-10 \\
39 & 18 & $(0.044, 0.405, 0.108, 0.333, 0.110)$ & $(5, 6, 4, 2)$ & 7.570e-15 & 1.317e-10 \\
40 & 18 & $(0.327, 0.057, 0.088, 0.386, 0.141)$ & $(7, 5, 3, 2)$ & 9.330e-14 & 8.192e-11 \\
\hline
\end{tabular}
\par\vspace{0.5em}
Errors larger than $10^{-6}$ are highlighted in bold.
\end{table}


\begin{table}[H]
\centering
\caption{Results from 40 randomized checks of the main theorem for $d=5$. Here $L$ denotes the left side and $R$ denotes the right side of the equation in Theorem~\ref{thm:normal.integral.representation}, and the relative error is $2|R-L|/(|R|+|L|)$.}
\label{tab:theorem-checks-d5}
\vspace{2em}
\begin{tabular}{r r l l r r}
\hline
Case & $n$ & $\bb{p}^+$ & $\bb{k}$ & Abs. error & Rel. error \\
\hline
1 & 34 & $(0.195, 0.214, 0.239, 0.033, 0.106, 0.213)$ & $(7, 3, 8, 8, 2)$ & 3.309e-09 & 2.740e-08 \\
2 & 25 & $(0.181, 0.265, 0.081, 0.191, 0.203, 0.079)$ & $(2, 2, 2, 6, 6)$ & 1.078e-07 & 1.131e-07 \\
3 & 22 & $(0.298, 0.188, 0.160, 0.095, 0.071, 0.188)$ & $(9, 3, 2, 2, 3)$ & 5.178e-10 & 5.672e-09 \\
4 & 14 & $(0.078, 0.059, 0.261, 0.395, 0.139, 0.069)$ & $(2, 2, 3, 3, 2)$ & 7.270e-10 & 1.192e-08 \\
5 & 28 & $(0.132, 0.154, 0.278, 0.068, 0.137, 0.230)$ & $(2, 6, 13, 2, 3)$ & 3.399e-11 & 5.701e-09 \\
6 & 17 & $(0.044, 0.189, 0.136, 0.180, 0.248, 0.203)$ & $(3, 2, 3, 2, 6)$ & 9.442e-11 & 2.029e-08 \\
7 & 28 & $(0.268, 0.129, 0.068, 0.048, 0.364, 0.123)$ & $(7, 10, 3, 5, 2)$ & 7.993e-13 & 4.436e-08 \\
8 & 30 & $(0.215, 0.150, 0.257, 0.140, 0.156, 0.082)$ & $(2, 4, 13, 4, 3)$ & 3.431e-08 & 8.052e-08 \\
9 & 13 & $(0.058, 0.081, 0.185, 0.164, 0.301, 0.211)$ & $(2, 3, 2, 2, 2)$ & 1.797e-11 & 2.880e-09 \\
10 & 32 & $(0.216, 0.065, 0.134, 0.079, 0.169, 0.338)$ & $(12, 8, 3, 3, 4)$ & 1.698e-15 & 1.620e-09 \\
11 & 27 & $(0.065, 0.259, 0.286, 0.067, 0.135, 0.189)$ & $(4, 5, 5, 3, 8)$ & 2.136e-10 & 1.803e-08 \\
12 & 20 & $(0.298, 0.123, 0.258, 0.110, 0.111, 0.100)$ & $(3, 3, 3, 2, 5)$ & 8.123e-09 & 9.418e-09 \\
13 & 22 & $(0.140, 0.100, 0.130, 0.432, 0.129, 0.069)$ & $(2, 7, 3, 4, 5)$ & 3.396e-11 & 1.494e-09 \\
14 & 16 & $(0.324, 0.227, 0.212, 0.034, 0.168, 0.036)$ & $(2, 3, 5, 2, 2)$ & 4.567e-08 & 5.939e-08 \\
15 & 22 & $(0.038, 0.391, 0.151, 0.067, 0.043, 0.310)$ & $(3, 4, 2, 7, 3)$ & 1.917e-11 & 4.032e-09 \\
16 & 18 & $(0.332, 0.110, 0.179, 0.116, 0.216, 0.046)$ & $(3, 3, 3, 2, 3)$ & 2.900e-08 & 3.595e-08 \\
17 & 15 & $(0.073, 0.143, 0.456, 0.051, 0.126, 0.150)$ & $(2, 2, 3, 3, 4)$ & 2.792e-10 & 3.288e-09 \\
18 & 22 & $(0.131, 0.180, 0.152, 0.412, 0.071, 0.053)$ & $(5, 3, 3, 3, 5)$ & 6.566e-10 & 6.533e-09 \\
19 & 26 & $(0.255, 0.169, 0.147, 0.219, 0.176, 0.034)$ & $(12, 3, 2, 2, 6)$ & 4.444e-11 & 4.066e-09 \\
20 & 23 & $(0.084, 0.109, 0.152, 0.129, 0.152, 0.375)$ & $(2, 7, 6, 2, 4)$ & 6.135e-14 & 5.597e-10 \\
21 & 15 & $(0.093, 0.156, 0.172, 0.048, 0.292, 0.239)$ & $(2, 2, 2, 3, 3)$ & 9.867e-11 & 9.526e-10 \\
22 & 13 & $(0.060, 0.224, 0.113, 0.314, 0.244, 0.044)$ & $(4, 2, 2, 2, 2)$ & 2.120e-11 & 9.828e-09 \\
23 & 19 & $(0.212, 0.109, 0.099, 0.112, 0.285, 0.182)$ & $(7, 2, 4, 2, 2)$ & 5.476e-11 & 1.329e-08 \\
24 & 24 & $(0.237, 0.108, 0.085, 0.255, 0.187, 0.128)$ & $(2, 2, 2, 11, 6)$ & 6.958e-09 & 5.101e-08 \\
25 & 19 & $(0.196, 0.331, 0.170, 0.196, 0.043, 0.064)$ & $(2, 2, 3, 8, 3)$ & 1.018e-08 & 1.691e-08 \\
26 & 23 & $(0.054, 0.140, 0.041, 0.183, 0.077, 0.505)$ & $(2, 6, 4, 2, 6)$ & 3.087e-13 & 9.298e-09 \\
27 & 16 & $(0.048, 0.150, 0.223, 0.295, 0.162, 0.124)$ & $(2, 3, 4, 3, 3)$ & 1.320e-11 & 5.857e-10 \\
28 & 12 & $(0.142, 0.095, 0.145, 0.259, 0.093, 0.266)$ & $(2, 2, 3, 2, 2)$ & 8.607e-11 & 2.671e-09 \\
29 & 22 & $(0.245, 0.362, 0.144, 0.124, 0.092, 0.033)$ & $(2, 7, 5, 2, 4)$ & 2.188e-09 & 2.483e-09 \\
30 & 16 & $(0.198, 0.155, 0.174, 0.202, 0.169, 0.102)$ & $(4, 4, 3, 2, 2)$ & 2.013e-10 & 4.250e-09 \\
31 & 34 & $(0.204, 0.293, 0.108, 0.152, 0.187, 0.055)$ & $(13, 2, 6, 3, 4)$ & 5.261e-10 & 4.489e-08 \\
32 & 21 & $(0.178, 0.238, 0.333, 0.053, 0.083, 0.115)$ & $(2, 6, 5, 2, 2)$ & 6.204e-10 & 1.012e-09 \\
33 & 15 & $(0.377, 0.170, 0.110, 0.104, 0.092, 0.148)$ & $(2, 3, 5, 2, 2)$ & 1.353e-09 & 5.478e-09 \\
34 & 19 & $(0.049, 0.271, 0.038, 0.254, 0.085, 0.302)$ & $(4, 2, 4, 5, 3)$ & 1.864e-12 & 7.881e-09 \\
35 & 13 & $(0.125, 0.165, 0.216, 0.079, 0.205, 0.210)$ & $(2, 3, 3, 2, 2)$ & 6.518e-10 & 1.518e-08 \\
36 & 33 & $(0.115, 0.095, 0.090, 0.113, 0.256, 0.331)$ & $(2, 10, 2, 3, 2)$ & 1.130e-10 & 5.570e-09 \\
37 & 32 & $(0.144, 0.088, 0.388, 0.146, 0.150, 0.084)$ & $(6, 10, 4, 4, 5)$ & 4.174e-11 & 5.535e-08 \\
38 & 20 & $(0.097, 0.063, 0.125, 0.175, 0.137, 0.403)$ & $(4, 2, 2, 3, 5)$ & 2.592e-11 & 3.506e-09 \\
39 & 29 & $(0.111, 0.064, 0.438, 0.037, 0.093, 0.257)$ & $(2, 3, 4, 9, 6)$ & 7.617e-09 & 5.288e-08 \\
40 & 21 & $(0.162, 0.212, 0.399, 0.039, 0.107, 0.081)$ & $(5, 4, 4, 2, 3)$ & 3.311e-09 & 2.001e-08 \\
\hline
\end{tabular}
\par\vspace{0.5em}
Errors larger than $10^{-6}$ are highlighted in bold.
\end{table}



\section*{Disclosure statement}
\addcontentsline{toc}{section}{Disclosure statement}

The author declares no conflicts of interest.

\section*{Funding}
\addcontentsline{toc}{section}{Funding}

F.\ Ouimet is supported by the Natural Sciences and Engineering Research Council of Canada through Discovery Grant RGPIN-2026-04471 and Discovery Launch Supplement DGECR-2026-00449.

\section*{References}
\addcontentsline{toc}{section}{References}

\setlength{\bibsep}{0pt plus 0ex}

\bibliographystyle{plainnat}
\bibliography{bib_clean}

\end{document}